\documentclass{amsart}

\usepackage{amscd}
\usepackage{amsmath}
\usepackage{amssymb}
\usepackage{amsthm}
\usepackage{epsf}
\usepackage{latexsym}
\usepackage{verbatim}

\input epsf.tex

\newtheorem{theorem}{Theorem}[section]
\newtheorem{definition}[theorem]{Definition}
\newtheorem{proposition}[theorem]{Proposition}
\newtheorem{lemma}[theorem]{Lemma}
\newtheorem{corollary}[theorem]{Corollary}

\theoremstyle{definition}

\newtheorem{example}[theorem]{Example}

\newcommand{\R}{\mathbb{R}}

\newcommand{\diam}{\operatorname{diam}}

\begin{document}

\title{Measurable Time-Restricted Sensitivity}

\author[Aiello]{Domenico Aiello}
\address[Domenico Aiello]{ Williams College, MA 01267, USA}
\email{Domenico Aiello <09da@williams.edu>}
\author[Diao] {Hansheng Diao}
\address[Hansheng Diao]{Massachusetts Institute of Technology, MA 02139, USA }
\email{hansheng@mit.edu}
\author[Fan]{Zhou Fan}
\address[Zhou Fan]{Harvard College, MA 02138, USA}
\email{Zhou Fan <zhoufan@fas.harvard.edu>}
\author[King]{Daniel O. King}
\address[Daniel King]{ Williams College, MA 01267, USA}
\email{Dan King <09dok@williams.edu>}
\author[Lin]{Jessica Lin}
\address[Jessica Lin]{Courant Institute of Mathematical Sciences, New York University, NY 10012, USA}
\email{jessicalin@nyu.edu}
\author[Silva]{Cesar E. Silva}
\address[Cesar Silva]{Department of Mathematics\\
     Williams College \\ Williamstown, MA 01267, USA}
\email{csilva@williams.edu}


\date{\today}

\maketitle 

\begin{abstract}
We develop two notions of time-restricted sensitivity to initial conditions for measurable dynamical systems, where the time before divergence of a pair of paths is at most an asymptotically logarithmic function of a measure of their initial distance. In the context of finite measure-preserving transformations on a compact space, we relate these notions to the metric entropy of the system. We examine one of these notions for classes of non-measure-preserving, nonsingular transformations.
\end{abstract}


\section{Introduction}\label{introduction}

Sensitivity has been widely studied as a characterization of chaos for topological dynamical systems, see e.g. \cite{BBCDS}, \cite{GW93}, \cite{A-B96}. Recently, sensitivity has been explored in the context of other measurable dynamical properties, such as weak mixing and entropy, for a finite measure-preserving transformation equipped with a metric of full support    \cite{ABC02},  \cite{He04}, \cite{Cadre05}. More recently, in \cite{J-S08}, the authors introduced a measure-theoretic version of sensitivity, invariant under measurable isomorphism, for nonsingular transformations. This has been further studied in \cite{G-S08}.

In broad terms, sensitivity asserts that for any point $x$ in the space, there exists another arbitrarily close point $y$ such that at some future positive time $n=n(x,y)$ the points $T^n(x)$ and $T^n(y)$ are separated by some predetermined distance. (Let us call this $n$ a {\it sensitive time}.) This simplest notion of sensitivity does not assert anything about the sensitive time other than its existence. A refinement of this definition, called strong sensitivity, was introduced in \cite{ABC02}, where for each point $x$ the set of sensitive times $n(x,y)$, for some $y$, is co-finite. Strong sensitivity was also studied in the measurable context in \cite{J-S08}.  An alternative refinement of (topological) sensitivity was studied in \cite{Moot07}, where the set of sensitive times is required to be syndetic. After completing this paper we learned
of \cite{Huang}, where the authors study extensions of pairwise sensitivity of for finite invariant measures
in topological dynamical systems.

In this paper, we are interested in placing a quantitative, asymptotic bound on the sensitive time, restricting the first sensitive time of a point $x$ and a point in an $\varepsilon$-ball around $x$ to be at most asymptotically logarithmic in the measure of that $\varepsilon$-ball. We develop two notions of measurable sensitivity that restrict the sensitive time in this way, show that these notions are related to positive metric entropy for finite measure-preserving systems, and explore one of the notions in the context of (non-measure-preserving) nonsingular transformations. An outstanding problem in nonsingular ergodic theory is the lack of a theory of entropy, see e.g. \cite{DS09}; our definitions of restricted sensitivity are related to entropy in the finite measure-preserving,  and can be
formulated in the context of nonsingular transformations and  be thought of as an
approach to positive nonsingular entropy.

In Section \ref{definitions}, we define the notions of restricted sensitivity, extending \cite{J-S08}, and restricted pairwise sensitivity, extending \cite{Cadre05}. We prove that under mild conditions restricted pairwise sensitivity implies restricted sensitivity, and we explore these notions for Bernoulli shifts. In Section \ref{entropy}, we consider the setting of finite measure-preserving transformations on a compact space and prove that restricted pairwise sensitivity implies positive metric entropy and positive metric entropy for a continuous, ergodic transformation implies restricted sensitivity. We explore a quantitative relationship between the entropy of a Bernoulli shift and the asymptotic bound on the sensitive time. In Section \ref{nonsingular}, we explore restricted sensitivity in the context of (non-measure-preserving) nonsingular transformations. We construct a class of nonsingular transformations, that includes  type III (i.e.,  not admitting an equivalent $\sigma$-finite invariant measure) transformations,  that are restricted sensitive as well as a class of nonsingular rank-one transformations that are not. It is well-known that
finite measure-preserving rank-one transformations have zero entropy; it would be interesting to know if
all nonsingular rank-one transformations are not restrictive sensitive. 

\subsection{Acknowledgements}
This paper is based on research by the Ergodic Theory group of the 2008 SMALL summer research project at Williams College.  Support for the project was provided by National Science Foundation REU Grant DMS - 0353634 and the Bronfman Science Center of Williams College.

\section{Time-Restricted Notions of Sensitivity}\label{definitions}

Throughout, let $(X,\mathcal{S}(X),\mu)$ denote a standard probability space. We also consider  a metric $d$ on $X$, and assume that $d$ is measurable so that the $d$-balls are measurable. It will be convenient for us to consider two compatibility conditions between $d$ and $\mu$. Let us say that $d$ is {\bf $\mu$-compatible} if all nonempty open $d$-balls have positive $\mu$-measure. The topology of a $\mu$-compatible metric is separable \cite{J-S08}. Let us say that $d$ is {\bf $\mu$-regular} if for all $x \in X$, there exists $c>0$ such that 
\[
c\mu(\overline{B_r}(x)) \leq \mu(B_r(x))\] for all $r>0$. In particular, if all $d$-balls are $\mu$-continuity sets, then $d$ is $\mu$-regular.

Recall the following definitions of measurable sensitivity:

\begin{definition} \cite{J-S08}\label{msdef}
A nonsingular dynamical system $(X,\mathcal{S}(X),\mu,T)$ is said to be
 \textbf{measurably sensitive} if whenever a dynamical
system $(X_1,\mathcal{S}(X_1),\mu_1,T_1)$ is measurably isomorphic
to $(X,\mathcal{S}(X),\mu,T)$ and $d$ is a $\mu_1$-supported metric on $X_1$, then
there exists $\delta > 0$ such that for all $x \in X_1$ and all  $\varepsilon > 0$, there exists
$n \in \mathbb{N}_0$ such that \[\mu_1 \{y\in B_{\varepsilon}(x):
d(T_1^nx,T_1^ny) > \delta \} > 0.\]
\end{definition}

\begin{definition} \cite{Cadre05}
A measure-preserving transformation $T$ on a probability space $(X,\mu)$ equipped with a metric $d$ is \textbf{pairwise sensitive} if there exists $\delta>0$ such that for $\mu^{\otimes^2}$-a.e. $(x,y) \in X \times X$, there exists $n \in \mathbb{N}_0$ such that \[d(T^nx,T^ny)>\delta.\] 
\end{definition}

In this paper, we study the following time-restricted modifications of these sensitivity notions:

\begin{definition}\label{rsdef}
A measurable, nonsingular transformation $T$ on the probability space $(X,\mathcal{S}(X),\mu)$ equipped with the metric $d$ is \textbf{restricted sensitive} if for a.e. $x \in X$, there exist $\delta>0$ and $a>0$ such that for every $\varepsilon>0$, there exists $n \in \mathbb{N}_0$, $n \leq -a \log \mu B_{\varepsilon}(x)$, with
\begin{equation*}
\mu\left\{y\in B_\varepsilon(x):d(T^nx,T^ny)>\delta\right\}>0.
\end{equation*}
\end{definition}

\begin{definition}\label{rpsdef}
A measurable, nonsingular transformation $T$ on the probability space $(X,\mathcal{S}(X),\mu)$ equipped with the metric $d$ is \textbf{restricted pairwise sensitive} if there exists $\delta>0$ and $a>0$ such that for $\mu^{\otimes^2}$-a.e. $(x,y) \in X \times X$, there exists $n\in\mathbb{N}_0$, $n \leq -a \log \mu B_{d(x,y)}(x)$ such that $d(T^nx,T^ny)>\delta$.
\end{definition}

Let us call $\delta$ and $a$ from these definitions a \textbf{sensitivity constant} and an \textbf{asymptotic rate}, respectively. For both Definitions \ref{rsdef} and \ref{rpsdef}, if $\delta$ is a sensitivity constant, then any $\delta'<\delta$ is also a sensitivity constant, and if $a$ is an asymptotic rate, than any $a'>a$ is also an asymptotic rate. Note that we allow the sensitivity constant and asymptotic rate to vary over $x \in X$ in Definition \ref{rsdef} but assume that they are constant across $x \in X$ in Definition \ref{rpsdef}. Finally, note that the condition $d(T^nx,T^ny)>\delta$ is true for $n=0$ when $d(x,y)>\delta$, so that it suffices to check $\varepsilon \leq \delta$ in Definition \ref{rsdef} and pairs of points $(x,y)$ with $d(x,y) \leq \delta$ in Definition \ref{rpsdef}. 

Restricted pairwise sensitivity is a stronger notion than restricted sensitivity, in the following sense:

\begin{proposition}\label{RPSimpliesRS}
Suppose $d$ is $\mu$-regular. If a transformation $T$ is restricted pairwise sensitive, then $T$ is restricted sensitive.
\end{proposition}

\begin{proof}
Let $\delta>0$ and $a>0$ be such that for a.e. $x \in X$, for a.e. $y \in X$, there exists $n \leq -a \log \mu B_{d(x,y)}(x)$ with $d(T^nx,T^ny)>\delta$. Consider any such $x \in X$, let $c>0$ be such that $c\mu(\overline{B_r}(x)) \leq \mu(B_r(x))$ for all $r>0$, and consider any $\varepsilon \leq \delta$. Let $\varepsilon_M=\sup\{\varepsilon'<\varepsilon:\mu(B_{\varepsilon'}(x))<\mu(B_\varepsilon(x))\}$. If $\mu(B_{\varepsilon_M}(x))=\mu(B_\varepsilon(x))$, then choose $\varepsilon'<\varepsilon_M$ sufficiently close to $\varepsilon_M$ such that $\mu(B_\varepsilon(x) \setminus B_{\varepsilon'}(x))>0$ and $\mu(B_{\varepsilon'}(x)) \geq c\mu(B_\varepsilon(x))$. Otherwise if $\mu(B_{\varepsilon_M}(x))<\mu(B_\varepsilon(x))$, then we must have $\mu(\overline{B_{\varepsilon_M}}(x))=\mu(B_\varepsilon(x))$, and we may choose $\varepsilon'=\varepsilon_M$ so that $\mu(B_\varepsilon(x) \setminus B_{\varepsilon'}(x))>0$ and $\mu(B_{\varepsilon'}(x)) \geq c\mu(\overline{B_{\varepsilon'}}(x))=\mu(B_\varepsilon(x))$ by the definition of $c$.

For a.e. $y \in B_\varepsilon(x) \setminus B_{\varepsilon'}(x)$, there exists $n \leq -a \log \mu B_{d(x,y)}(x)$ such that $d(T^nx,T^ny)>\delta$. We note that $\mu B_{d(x,y)}(x) \geq \mu B_{\varepsilon'}(x) \geq c \mu B_\varepsilon(x)$, so $n \leq -a \log c \mu B_\varepsilon(x)=-a \log \mu B_\varepsilon(x)-a \log c$. As we may choose the pairwise sensitivity constant $\delta$ so that $\mu B_\varepsilon(x)\leq\mu B_\delta(x)<1$, we may take $\tilde{a}$ so that $-a \log \mu B_\varepsilon(x)-a \log c \leq -\tilde{a} \log \mu B_\varepsilon(x)$ for all $\varepsilon \leq \delta$. Hence $T$ is restricted sensitive with sensitivity constant $\delta$ and sensitivity function $\tilde{a}$.
\end{proof}

As an application of these notions of restricted sensitivity and restricted pairwise sensitivity, let us consider the standard one-sided and two-sided Bernoulli shift transformations.

\begin{example}\label{onesidedshift}
Consider the space $\Sigma_N^+=\prod_{i=0}^\infty \{1,\ldots,N\}$ with its product $\sigma$-algebra and the probability measure $\mu=\bigotimes_{i=0}^\infty \mu'$, where $\mu'$ is a probability measure on $\{1,\ldots,N\}$ of full support. Consider the metric $d(\sigma,\tau)=2^{-I(\sigma,\tau)}$ where $I(\sigma,\tau)=\min \{i \geq 0:\sigma_i \neq \tau_i\}$. Let $T$ be the one-sided Bernoulli shift transformation $T((\sigma_0,\sigma_1,\sigma_2,\ldots))=(\sigma_1,\sigma_2,\sigma_3,\ldots)$. Let $p=\max_k \mu'(k)$, $\delta=\frac{1}{2}$, and $a=-\frac{1}{\log p}$. Consider any two points $\sigma,\tau \in \Sigma_N^+$ and suppose that $I(\sigma,\tau)=n$. Then $\mu B_{d(\sigma,\tau)}(\sigma) \leq p^n$, so $-a\log \mu B_{d(\sigma,\tau)}(\sigma) \geq n$. Since $d(T^n\sigma,T^n\tau)=1>\delta$, $T$ is restricted pairwise sensitive. As $d$ is $\mu$-regular with $c=\min_k \mu'(k)$ for all $\sigma \in \Sigma_N^+$, $T$ is also restricted sensitive.
\end{example}

\begin{example}\label{twosidedshift}
Consider the space $\Sigma_N=\prod_{i=-\infty}^\infty \{1,\ldots,N\}$ with its product $\sigma$-algebra and the probability measure $\mu=\bigotimes_{i=-\infty}^\infty \mu'$, where $\mu'$ is a probability measure on $\{1,\ldots,N\}$ of full support. Consider the metric $d(\sigma,\tau)=2^{-I(\sigma,\tau)}$, where $I(\sigma,\tau)=\min \{|i|: \sigma_i \neq \tau_i\}$. Let $T$ be the two-sided Bernoulli shift transformation such that $T(\sigma)_i=\sigma_{i+1}$. Let $p_k=\mu'(k)$ for each $k=1,\ldots,N$, and consider any $\delta<1$ and any $a>0$. For any $\sigma \in \Sigma_N$, choose an integer $k_1>0$ such that $2^{-k_1}<\delta$ and let $P=\prod_{i=-k_1+1}^{k_1-1} p_{\sigma_i}$. Choose an integer $k_2>k_1$ such that $2^{-(k_2+a \log P)}<\delta$. Consider any $\tau$ in the cylinder set $[\bar{\sigma}_{-k_1}\sigma_{-k_1+1}\ldots \sigma_0 \ldots \sigma_{k_2}]$, where $\bar{\sigma}_{-k_1}$ is some symbol not equal to $\sigma_{-k_1}$. Then by the construction of $k_1$ and $k_2$, for all $n \leq -a \log P=-a \log \mu B_{d(\sigma,\tau)}(\sigma)$, $d(T^n\sigma,T^n\tau)\leq \max(2^{-k_1},2^{-(k_2-n)}) \leq \max(2^{-k_1},2^{-(k_2+a \log P)})<\delta$. Hence $T$ is not restricted pairwise sensitive.

On the other hand, let $p=\max_k \sigma'(k)$, $\delta=\frac{1}{4}$, and $a=-\frac{1}{2\log p}$. Consider any $\sigma \in \Sigma_N$ and any ball $B(\sigma)=[\sigma_{-n}\ldots \sigma_0 \ldots \sigma_n]$. We have that $\mu B(\sigma) \leq p^{2n+1}$, so $n<-a \log \mu B(\sigma)$. Since $d(T^n\sigma,T^n\tau)=\frac{1}{2}>\delta$ for all $\tau \in [\sigma_{-n}\ldots \sigma_n\bar{\sigma}_{n+1}] \subset [\sigma_{-n}\ldots \sigma_n]$ where $\bar{\sigma}_{n+1}$ is some symbol not equal to $\sigma_{n+1}$, $T$ is restricted sensitive.
\end{example}

Thus we have examples of transformations that are restricted pairwise sensitive and transformations that are restricted sensitive but not restricted pairwise sensitive. One can also easily construct transformations that are not restricted sensitive but that are measurably sensitive according to Definition \ref{msdef}; examples of such transformations include the rank-one cutting and stacking transformations that we will study in Section \ref{nonsingular}.

\section{Sensitivity and Entropy for Measure-Preserving Transformations}\label{entropy}

In the setting of finite measure-preserving transformations on a compact metric space, the notion of metric entropy measures the rate at which a transformation disorganizes the space. As measurable time-restricted notions of sensitivity convey that points separate from one another rapidly on a local level, it is natural to explore the connection between notions of measurable time-restricted sensitivity and the metric entropy of a dynamical system.

The following theorem and its immediate corollary show that, under mild conditions on $\mu$ and $d$, restricted pairwise sensitivity implies positive metric entropy in the context of finite measure-preserving transformations on a compact metric space:

\begin{theorem} \label{RPSimpliesPosEntropy}
Suppose $\mu$ is a nonatomic probability measure on $X$ and $d$ is $\mu$-compatible and $\mu$-regular. Let $T$ be a measure-preserving transformation on $X$. Let $\mathcal{A}=\{A_1,\ldots,A_k\}$ be a partition of $X$ such that $\diam A_i<\delta$ for all $i$, and let $h_\mu(T,\mathcal{A})$ be the metric entropy of $T$ with respect to $\mathcal{A}$. If $T$ is restricted pairwise sensitive with asymptotic rate $a$, then $h_\mu(T,\mathcal{A}) \geq \frac{1}{a}$.
\end{theorem}

\begin{proof}
Suppose $T$ is restricted pairwise sensitive with sensitivity constant $\delta$ and asymptotic rate $a$. Take $x \in X$ such that for a.e. $y \in X$, there exists $n \leq -a\log \mu B_{d(x,y)}(x)$ with $d(T^nx,T^ny)>\delta$. Let $C_n(x)$ denote the element of the partition $\bigvee_{i=0}^n T^{-i}\mathcal{A}$ containing $x$. If $y \in C_n(x)$, then $T^ix$ and $T^iy$ are in the same element of $\mathcal{A}$ for all $i \leq n$, so $d(T^ix,T^iy)<\delta$ for all $i \leq n$.

Suppose $c>0$ is such that $e^{-c}\mu(\overline{B_r}(x)) \leq \mu(B_r(x))$ for all $r>0$. (The existence of this $c$ is given by $\mu$-regularity.) Since $\mu$ is nonatomic, $\mu(B_r(x)) \to 0$ as $r \to 0$. For each integer $n \geq 1$, we may take $\varepsilon_n$ such that $e^{-(n-1)c}>\mu(B_{\varepsilon_n}(x)) \geq e^{-nc}$, for otherwise there must exist $\varepsilon'$ with $\mu(\overline{B_{\varepsilon'}}(x)) \geq e^{-(n-1)c}$ and $\mu(B_{\varepsilon'}(x))<e^{-nc}$, contradicting the definition of $c$. If $y \notin B_{\varepsilon_{\lfloor \frac{n}{ac} \rfloor}}(x)$, then $\mu B_{d(x,y)}(x) \geq B_{\varepsilon_{\lfloor \frac{n}{ac} \rfloor}}(x) \geq e^{-\lfloor \frac{n}{ac} \rfloor c}$, so for a.e. $y \notin B_{\varepsilon_{\lfloor \frac{n}{ac} \rfloor}}(x)$, there exists $i \leq -a\log \mu B_{d(x,y)}(x) \leq n$ such that $d(T^ix,T^iy)>\delta$, and so $y \notin C_n(x)$. Hence $C_n(x) \subset B_{\varepsilon_{\lfloor \frac{n}{ac} \rfloor}}(x) \bmod \mu$, so $\mu C_n(x) \leq \mu B_{\varepsilon_{\lfloor \frac{n}{ac} \rfloor}}(x)<e^{(\lfloor \frac{n}{ac} \rfloor-1)c}$. Letting $h(x)=\liminf_{n \to \infty} -\frac{1}{n+1} \log \mu C_n(x)$, we have
\[h(x) \geq \liminf_{n \to \infty} -\frac{c}{n+1}\left(\left\lfloor \frac{n}{ac} \right\rfloor-1\right)=\frac{1}{a}.\]
This holds for a.e. $x \in X$, so by Fatou's lemma,
\[h_\mu(T,\mathcal{A})=\lim_{n \to \infty} \int -\frac{1}{n+1} \log \mu C_n(x)\,d\mu \geq \int h(x)\,d\mu \geq \frac{1}{a}.\]
\end{proof}

\begin{corollary}
Suppose $(X,d)$ is compact, $\mu$ is nonatomic, and $d$ is $\mu$-compatible and $\mu$-regular. Let $T$ be a measure-preserving transformation on $X$. If $T$ is restricted pairwise sensitive with asymptotic rate $a$, then $h_\mu(T) \geq \frac{1}{a}$.
\end{corollary}

Example \ref{twosidedshift} shows that a converse to this theorem is not true, for the two-sided Bernoulli shift has positive entropy but is not restricted pairwise sensitive. The following result shows, however, that the implication in this direction is true if we replace restricted pairwise sensitivity with restricted sensitivity, under assumptions of ergodicity and continuity of $T$:

\begin{theorem}\label{PosEntropyImpliesRS}
Suppose $(X,d)$ is compact and $\mu$ is nonatomic. Let $T$ be a continuous, ergodic, measure-preserving transformation on $X$. If $h_\mu(T)>\frac{1}{a}$, then $T$ is restricted sensitive with asymptotic rate $a$ for all $x \in X$.
\end{theorem}

\begin{proof}
Suppose $T$ is not restricted sensitive with asymptotic rate $a$ over all $x \in X$. Then there exists a positive measure set $A$ such that for any $x \in A$ and $\delta>0$, there exists $\varepsilon(\delta) \leq \delta$ such that $\mu\left\{y\in B_{\varepsilon(\delta)}(x): d(T^nx,T^ny)>\delta\right\}=0$ for all $n \leq -a \log \mu B_{\varepsilon(\delta)}(x)$. For each $n$, let
\[C(x,n,\delta)=\{y\in X: d(T^i x,T^i y) \leq \delta \text{ for all } 0\leq i<n\}.\]
Then $B_{\varepsilon(\delta)}(x)\subset C(x,\left\lceil -a \log \mu B_{\varepsilon(\delta)}(x)\right\rceil,\delta) \mod \mu$, so
\[\mu B_{\varepsilon(\delta)}(x) \leq \mu C(x,\left\lceil -a\log\mu B_{\varepsilon(\delta)}(x)\right\rceil,\delta).\]
As $\mu$ is nonatomic, $\lim_{\delta \to \infty} -a \log \mu B_{\varepsilon(\delta)}(x)=\infty$. Then using the Brin-Katok Theorem, there exists an $x \in A$ such that the metric entropy is given by
\begin{align*}
h_\mu(T)&=\lim_{\delta\rightarrow 0} \liminf_{n\rightarrow\infty} -\frac{1}{n} \log \mu C(x,n,\delta)\\
&\leq \lim_{m \rightarrow \infty} \liminf_{n \rightarrow \infty} \frac{-\log \mu C(x,\lceil -a \log \mu B_{\varepsilon(\frac{1}{n})}(x) \rceil,\frac{1}{m})}{\lceil -a \log \mu B_{\varepsilon(\frac{1}{n})}(x) \rceil}\\
&\leq \liminf_{n \rightarrow \infty} \frac{-\log \mu C(x,\lceil -a \log \mu B_{\varepsilon(\frac{1}{n})}(x) \rceil,\frac{1}{n})}{\lceil -a \log \mu B_{\varepsilon(\frac{1}{n})}(x) \rceil}\\
&\leq \liminf_{n \rightarrow \infty} \frac{-\log \mu B_{\varepsilon(\frac{1}{n})}(x)}{\lceil -a \log \mu B_{\varepsilon(\frac{1}{n})}(x) \rceil}\\
&\leq \frac{1}{a},
\end{align*}
where in the third line we have used the inequality $\mu C(x,n,\delta) \leq \mu C(x,n,\delta')$ if $\delta \leq \delta'$. This gives the desired contradiction, and thus $T$ is restricted sensitive with asymptotic rate $a$.
\end{proof}

The conditions of ergodicity and continuity on $T$ in the above proof were required for the use of the Brin-Katok Theorem. Note that ergodicity of $T$ is necessary for any result of this type, because positive entropy is a global property of the system whereas sensitivity is a condition that must hold locally at almost every point. (Indeed, if we consider the disjoint union of a restricted sensitive system with positive entropy and a non-restricted sensitive system with zero entropy, the resulting system would have positive entropy but not be restricted sensitive.) That the continuity of $T$ is necessary is not clear, and we suspect that a result of this type holds true without the continuity condition.

Regarding the converse of this theorem, the following example constructs an ergodic, measure-preserving transformation which is restricted sensitive but has zero measure-theoretic entropy:

\begin{example}
Let $X$ be a disjoint union of two copies of $[0,1/2]$, labeled $I_1$ and $I_2$, equipped with their Borel sigma algebras $\mathcal{B}_1$ and $\mathcal{B}_2$ and Lebesgue measures. Consider the $\sigma$-algebra on $X$ given by $\mathcal{S}=\{S_1 \cup S_2 : S_1 \in \mathcal{B}_1 \textrm{ and } S_2 \in \mathcal{B}_2\}$ and probability measure $\mu$ on $X$ given by $\mu(A) = \lambda(A \cap I_1) + \lambda(A \cap I_2)$. Define a metric $d$ on $X$ by $d(x,y) = |x-y|$ if $x$ and $y$ are in the same copy of $[0, 1/2]$, and $d(x,y) = 2$ if not.

Let $A \subset [0,1/2]$ be a Borel set of Lebesgue measure $\lambda(A)=\frac{1}{4}$ with the property that for any subinterval $K$ of $[0,1/2]$, $\lambda(K\cap A)>0$ and $\lambda(K \cap A^c)>0$. (For a construction of such a set, see Appendix B of \cite{Si08}.) By Theorem 3.4.23 in \cite{srivastava1998cbs}, there exist measure-preserving Borel isomorphisms $\phi: A \rightarrow [0,1/4)$ and $\psi: A^c \rightarrow [0,1/4)$. Let $A_1$ and $A_2$ be the copies of $A$ inside $I_1$ and $I_2$, respectively, and let $\phi_1,\psi_1$ and $\phi_2,\psi_2$ be copies of the maps $\phi, \psi$ on $I_1,I_2$. Define a transformation $T : X \rightarrow X$ by
\[
T(x) = \left\{
 \begin{array}{ll}
 \psi_1^{-1} \circ \phi_1(x) & \textrm{ for } x \in A_1\\
 \phi_2^{-1} \circ \psi_1(x) & \textrm{ for } x \in A_1^c \cap I_1\\
 \psi_2^{-1} \circ \phi_2(x) & \textrm{ for } x \in A_2\\
 \phi_1^{-1} \circ R \circ \psi_2(x) & \textrm{ for } x \in A_2^c \cap I_2,
 \end{array}
\right.
\]
where $R: [0,1/4) \rightarrow [0,1/4)$ is an irrational rotation.  Note that $T$ is finite measure-preserving.  Also, note that when $T^4$ is restricted to any one of the four ``segments'' $A_1$, $A_1^c \cap I_1$, $A_2$, or $A_2^c \cap I_2$, it is isomorphic to an irrational rotation. We claim that $T$ is ergodic and restricted sensitive and that $h_\mu(T) = 0$.

To see that $T$ is ergodic, let $C \subset X$ and $D \subset X$ have positive measure.  Then we can find positive measure subsets $C^*\subset C$ and $D^* \subset D$ each of which is completely contained inside one of the four segments $A_1$, $A_1^c \cap I_1$, $A_2$, and $A_2^c \cap I_2$.  Let $k < 4$ be the integer such that $T^k(C^*)$ and $D^*$ are in the same segment.  Then, since $T^4$ is isomorphic to an irrational rotation when restricted to that segment, there exists $n$ a multiple of 4 such that $\mu(T^{k+n}(C^*) \cap D^*) > 0$. Therefore, $\mu(T^{k+n}(C) \cap D) > 0$.

To see that $T$ is restricted  sensitive, choose $\delta = 1$ and $a = 2$ for all $x \in X$. For any $x \in I_1$ and $\varepsilon \leq 1$, by construction, $\mu(B_\varepsilon(x) \cap A_1)>0$ and $\mu(B_\varepsilon(x) \cap (A_1^c \cap I_1))>0$, so $\mu\{y \in B_\varepsilon(x):d(Tx,Ty)>\delta\}>0$ and we note that $n=1<-2\log \frac{1}{2} \leq -a\log\mu B_\varepsilon(x)$. A similar argument holds for any $x \in I_2$, so $T$ is restricted sensitive.

To see that $T$ has zero entropy, simply note that $h_\mu(T^4)=0$ because $T^4$ is isomorphic to a disjoint union of 4 irrational rotations. This implies that $h_\mu(T)=0$.
\end{example}

Theorem \ref{RPSimpliesPosEntropy} places a quantitative lower bound on the entropy of a system using the asymptotic rate $a$ in Definition \ref{rpsdef}. In the context of restricted sensitivity, we may consider the minimal asymptotic rate $a_T^*(x)=\inf_{\delta(x)} \inf_{a(x)} a(x)$ of a restricted sensitive transformation $T$ as a function of $x \in X$, where the infimums are taken over all sensitivity constants and asymptotic rates at the point $x$. Theorem \ref{PosEntropyImpliesRS} then implies that if $a>0$ is such that $a_T^*(x)>a$ for a.e. $x \in X$, then $\frac{1}{a} \geq h_\mu(T)$. This function $a_T^*(x)$ is well-defined over a.e. $x \in X$, and it is in fact measurable under mild conditions on $\mu$ and $d$:

\begin{proposition}\label{ameasurable}
Suppose $(X,d)$ is separable and $\mu$ is a Borel probability measure. Then for any restricted sensitive transformation $T$ on $X$, its minimal asymptotic rate function $a_T^*:X \to \R$ is measurable.
\end{proposition}

We defer the proof of this technical proposition to Appendix \ref{appendix}. The following proposition computes $a_T^*$ for the one-sided Bernoulli shift transformation of Example \ref{onesidedshift} and shows that the upper bound on the entropy of the Bernoulli shift from Theorem \ref{PosEntropyImpliesRS} is tight:

\begin{proposition}
Let $(\Sigma_N^+,\mathcal{B},\mu)$, $d$, and $T$ be as in Example \ref{onesidedshift}. Then \[\frac{1}{a_T^*(\sigma)}=h_\mu(T)\] for a.e. $\sigma \in \Sigma_N^+$.
\end{proposition}

\begin{proof}
Take $\sigma \in \Sigma_N^+$ for which the restricted sensitivity condition holds with sensitivity constant $\delta(\sigma)>0$ and asymptotic rate $a(\sigma)>0$. Let $c(\sigma) \geq 0$ be the integer such that $2^{-c(\sigma)}>\delta(\sigma) \geq 2^{-c(\sigma)-1}$. For each $i=1,\ldots,N$, let $k_i^{(n)}(\sigma)$ be the number of occurrences of the symbol $i$ in $\sigma_0$ through $\sigma_{n-1}$. For any ball $B(\sigma)=[\sigma_0,\sigma_1,\sigma_2, \ldots \sigma_{n-1}]$, $\min \{k:\mu(\tau \in B(\sigma):d(T^k\sigma,T^k\tau)>\delta(\sigma))>0\}=n-c(\sigma)$. Hence we must have that $n-c(\sigma) \leq -a(\sigma) \log(p_1^{k_1^{(n)}(\sigma)}\ldots p_N^{k_N^{(n)}(\sigma)})$ for all $n>c(\sigma)$, so for fixed $\delta(\sigma)$,
\[\inf a(\sigma)=\sup_{n>c(\sigma)} -(n-c(\sigma))\left(\log(p_1^{k_1^{(n)}(\sigma)}\ldots p_N^{k_N^{(n)}(\sigma)})\right)^{-1}.\]
Thus
\begin{align*}
a_T^*(\sigma)&=\inf_{\delta(\sigma)} \sup_{n>c(\sigma)} -(n-c(\sigma))\left(\log(p_1^{k_1^{(n)}(\sigma)}\ldots p_N^{k_N^{(n)}(\sigma)})\right)^{-1}\\
&=\left(\sup_{c \geq 0} \inf_{n>c} \frac{n}{n-c} \cdot \frac{1}{n}\sum_{i=1}^N -k_i^{(n)}(\sigma)\log p_i\right)^{-1}\\
&=\left(\lim_{c \to \infty} \inf_{n>c} \frac{n}{n-c} \cdot \frac{1}{n}\sum_{i=1}^N -k_i^{(n)}(\sigma)\log p_i\right)^{-1}.
\end{align*}
By the Birkhoff Ergodic Theorem, for a.e. $\sigma \in \Sigma_N^+$,
\[\lim_{n \to \infty} \frac{k_i^{(n)}(\sigma)}{n}=\lim_{n \to \infty}\frac{1}{n}\sum_{i=1}^N \chi_{[i]}\left(T^i(\sigma)\right)p_i,\]
so
\[\lim_{n \to \infty} \frac{1}{n}\sum_{i=1}^N -k_i^{(n)}(\sigma)\log p_i=\sum_{i=1}^N -p_i\log p_i=h_\mu(T).\]
Then
\[\frac{1}{a_T^*(\sigma)} \leq \lim_{c \to \infty} \lim_{n \to \infty} \frac{n}{n-c} \cdot \frac{1}{n}\sum_{i=1}^N -k_i^{(n)}(\sigma)\log p_i=h_\mu(T).\]
For any $\varepsilon>0$, there exists $c$ sufficiently large so that $\frac{1}{n}\sum_{i=1}^N -k_i^{(n)}(\sigma)\log p_i \geq h_\mu(T)-\varepsilon$ for all $n>c$; hence
\[\frac{1}{a_T^*(\sigma)} \geq \lim_{c \to \infty} \lim_{n \to \infty} \frac{n}{n-c} \cdot (h_\mu(T)-\varepsilon)=h_\mu(T)-\varepsilon.\]
As $\varepsilon$ was arbitrary, $\frac{1}{a_T^*(\sigma)}=h_\mu(T)$.
\end{proof}

\section{Restricted Sensitivity for Nonsingular Transformations}\label{nonsingular}

In the preceding section, the notions of restricted sensitivity and restricted pairwise sensitivity were used to examine measure-preserving transformations, for which there is a well-developed theory of metric entropy. These notions can be applied as well to nonsingular transformations; let us consider restricted sensitivity in this section. We show that a general class of nonsingular rank-one transformations (including measure-preserving rank-one transformations) are not restricted sensitive, and we construct a class of nonsingular type III transformations that are restricted sensitive.

Let us recall the definition of rank-one transformations. This class is known to contain finite measure-preserving mixing transformations \cite{Or70} and type III power weakly mixing nonsingular transformations \cite{AFS01}.  A nonsingular transformation $T: (X,\mu) \rightarrow (X,\mu)$
is \textbf{type III} if there are no $\sigma$-finite measures invariant under $T$ that are equivalent to $\mu$. The first example of a type III transformation was rank-one \cite{Or60}. By~\cite{J-S08}, it follows that the class of rank-one transformations further includes strong measurably sensitive finite measure-preserving transformations and measurably sensitive type III transformations. 

We give the cutting and stacking definition of these transformations and follow the notation of \cite{CS04} and \cite{DS09}. Our presentation includes nonsingular transformations. A \textbf{column} consists of a finite ordered collection of disjoint intervals in $\mathbb{R}$.  Each interval is called a \textbf{level}, and the levels may be of different lengths. The \textbf{height} of the column is the number of levels in the column.  Each column defines an  associated \textbf{column map}, defined on all levels except the top, by mapping each interval of the column to the next  interval in the column by the unique orientation-preserving affine map that takes one interval to the other. Hence the column map is defined on all but the last  level.

A rank-one nonsingular transformation is specified by a sequence of integers 
  $\{r_n \geq 2\}$,    a sequence of functions $s_n:\{0,\ldots, r_n-1\}\to \mathbb N_0$,  
  and a sequence  $\{p_{n}\}$  of probability vectors on $\{0,\ldots, r_n-1\}$.  
  In the case of a measure-preserving transformations the probability vectors 
  are all uniform, i.e.,  $p_n(i)=1/r_n$ for all $i\in \{0,\ldots, r_n-1\}$. 
  We now describe the inductive  procedure that constructs a sequence of columns $C_{n}$. Start by letting $C_0$ consist of a single interval. 
Assume that column   $C_{n} = \{I_{n,i}\}_{i=0}^{h_{n}-1}$ of  height $h_{n}$ has been constructed.
This  we have a  column
map $T_{n}$,
where $T_{n}(I_{n,i}) = I_{n,i+1}$
for $i\ne h_{n}-1$. To construct column $C_{n+1}$ subdivide     $C_{n}$ into $r_{n}$ subcolumns by
cutting each level $I_{n,i}$ into $r_{n}$ subintervals or  \textbf{sublevels},
$\{I_{n,i}^{[j]}\}_{j=0}^{r_{n}-1}$,
(where
$I_{n,i}^{[0]}$ is the leftmost sublevel and $I_{n,i}^{[r_{n}-1]}$ is
the rightmost) whose lengths are in the proportions \[p_n(0) : p_n(1) : \cdots : p_n(r_n-1).\]
(For example, if $r_n=1$, and $p_n(0)=1/3, p_n(1)=2/3$, then every level is cut in the 
proportions $1/3 : 2/3$.)
Then  the \textbf{subcolumns} of $C_{n}$ are $C_{n}^{[j]} =
\{I_{n,i}^{[j]}\}_{i=0}^{h_{n}-1}$.  By
preserving the order on the levels, each subcolumn, $C_{n}^{[j]}$, is
a column in its own right with the associated map $T_{n}^{[j]}$ which is the
restriction of $T_{n}$ to $C_{n}^{[j]}$.
The next step is to place new intervals the size of the top sublevel 
 above each subcolumn
by adding $s_{n,j}$ levels above $C_{n}^{[j]}$; these new intervals are called
 \textbf{spacer levels}. To obtain the next column then we 
\textbf{stack} the resulting subcolumns with spacers right on top of left
yielding the  new column $C_{n+1}$ with height \[h_{n+1} = r_{n}h_{n} +
\sum_{j=0}^{r_{n}-1}s_{n}(j).\]
 Let $S_{n}$  denote the \textbf{union of spacer
levels} added to $C_{n}$, the collection of levels in
$C_{n+1}$ that are not sublevels of levels in $C_{n}$; hence,
$C_{n+1} = C_{n} \sqcup S_{n}$, and denote by $S_{n}^{[j]}$ the collection of spacers
added over $C_{n}^{[j]}$, so that $S_n=\bigcup_{j=0}^{r_n-1} S_{n}^{[j]}$.     The associated column map $T_{n+1}$
restricts to $T_{n}$ on the levels in $C_{n}$.   Let $X$  be the union of all the levels in all columns. We assume that as $n\to\infty$ the maximal length of the intervals in $C_n$
converges to $0$, so we may define a transformation $T$ of $(X,\mu)$ by 
\[T(x):=\lim_{n\to\infty} T_{{n}}(x).\]
One can verify that $T$ is well-defined and invertible  a.e. and that it is  nonsingular and ergodic.  $T$ is measure-preserving if all the probability vectors $p_n$ are uniform, and $\mu(X)<\infty$ if and only if the total measure of the added spacers is finite.

The following proposition shows that nonsingular rank-one transformations constructed in this way are not restricted sensitive if there is a uniform lower bound on the elements of the probability vectors $\{p_n\}$.

\begin{proposition}\label{notrestsens}
Let $T$ be a nonsingular rank-one transformation on $[0,1)$ with the Euclidean metric $d$ and Lebesgue measure $\lambda$. Suppose that column $C_n$ is divided into $r_n$ subcolumns with proportions $p_n(0),\ldots,p_n(r_n-1)$. If there exists $c>0$ such that $p_n(j) \geq c$ for all $n$ and $j$, then $T$ is not restricted sensitive.
\end{proposition}

\begin{proof}
For each $n$, let us further divide the leftmost subcolumn $C_n^{[0]}$ of column $C_n$ into three equal subcolumns, labeled from left to right as $C_n^{[0],1}$, $C_n^{[0],2}$, and $C_n^{[0],3}$. Let $S_n=C_n \setminus C_n^{[0],2}$ and let $S=\bigcup_{k=0}^\infty \bigcap_{n=k}^\infty S_n$. As $\lambda S_n \leq (1-\frac{c}{3})\lambda C_n \leq 1-\frac{c}{3}$, $\lambda(\bigcap_{n=k}^\infty S_n) \leq 1-\frac{c}{3}$ for each $k$ and thus $\lambda S \leq 1-\frac{c}{3}$. So the complement of $S$ in $[0,1)$ has positive measure.

For any $x \notin S$, there is an increasing sequence $\{n_k\}$ such that $x \in C_{n_k} \setminus S_{n_k}$. Let $h_n$ and $w_n$ be the height and the width of the smallest level, respectively, of column $C_n$. For any $\delta>0$ and $a>0$, there exists $n_k$ sufficiently large such that $w_{n_k}<\delta$ and $a(n_k \log \frac{1}{c}+\log \frac{3}{2w_0})<2^{n_k-1}$. By the construction of $S_{n_k-1}$, if $h$ is the smallest number such that $T^h(x)$ is in the highest level of $C_{n_k}$ and $w$ is the distance from $x$ to the closer of the two endpoints of the level containing $x$ in $C_{n_k}$, then $h \geq \frac{h_{n_k}}{2}$ and $w \geq \frac{w_{n_k}}{3}$. We also note that $w_{n_k} \geq c^{n_k}w_0$ and $h_{n_k} \geq 2^{n_k}$.

Consider the ball $B_w(x)$. We note that
\begin{align*}
-a \log \lambda B_w(x)&=-a \log 2w \leq -a \log \frac{2w_{n_k}}{3} \leq -a \log \frac{2c^{n_k}w_0}{3}\\
&=a\left(n_k \log \tfrac{1}{c}+\log \tfrac{3}{2w_0}\right)<2^{n_k-1}\leq \frac{h_{n_k}}{2} \leq h.
\end{align*}
Hence, for all $y \in B_w(x)$, $T^ny \in B_w(T^nx)$ for all $n \leq -a\log \lambda B_w(x)$, so \[d(T^nx,T^ny)<w<\delta.\] Hence $T$ is not restricted  sensitive.
\end{proof}

This proposition addresses a large class of nonsingular rank-one transformations. As measure-preserving transformations are those for which the probability vectors $p_n$ are uniform, measure-preserving transformations for which $r_n$ is bounded above over all $n$ satisfy the conditions of this proposition. In fact, the argument in the above proof can be modified to hold for all measure-preserving rank-one transformations:

\begin{proposition}
If $T$ is a measure-preserving rank-one transformation on $[0,1)$ with the Euclidean metric $d$ and Lebesgue measure $\lambda$, then $T$ is not restricted sensitive.
\end{proposition}

\begin{proof}
As the proof is very similar to that for Proposition \ref{notrestsens}, we will highlight the modifications required. Let us divide each subcolumn $C_n^{[j]}$ of column $C_n$ into three equal subcolumns, labeled from left to right as $C_n^{[j],1}$, $C_n^{[j],2}$, and $C_n^{[j],3}$. Let
\[S_n=\bigcup_{j=\lceil \frac{r_n-1}{2} \rceil}^{r_n-1} C_n^{[j]} \cup \bigcup_{j=0}^{\lceil \frac{r_n-1}{2}\rceil-1} C_n^{[j],1} \cup \bigcup_{j=0}^{\lceil \frac{r_n-1}{2}\rceil-1} C_n^{[j],3},\]
and let $S=\bigcup_{k=0}^\infty \bigcap_{n=k}^\infty S_n$ as before. We have $\lambda S_n \leq \frac{8}{9} \lambda C_n \leq \frac{8}{9}$ for all $n$, so the complement of $S$ in $[0,1)$ has positive measure.\\

The rest of the proof is the same as for Proposition \ref{notrestsens}, except that for any $\delta>0$ and $a>0$, we choose $n_k$ sufficiently large such that $w_{n_k}<\delta$ and $a \log \frac{3h_{n_k}}{2w_0h_0}<\frac{h_{n_k}}{2}$. As $w_{n_k}h_{n_k} \geq w_0h_0$ for all $k$, we use the bound $w_{n_k} \geq \frac{w_0h_0}{h_{n_k}}$ in place of the bounds $w_{n_k} \geq c^{n_k}w_0$ and $h_{n_k} \geq 2^{n_k}$ from the proof of Proposition \ref{notrestsens}.
\end{proof}

We end by constructing type III nonsingular transformations that are restricted sensitive. The construction is a consequence of the following proposition:

\begin{proposition}\label{prod}
Let $T$ be a transformation on a probability space $(X,\mu)$ with metric $d_X$ and $S$ be a transformation on a probability space $(Y,\nu)$ with metric $d_Y$, and let $d_Y$ be $\nu$-supported. If $T$ is restricted sensitive, then the transformation $T \times S$ on $X \times Y$ (with the product $\sigma$-algebra and product measure) is restricted sensitive under the metric given by
\begin{equation*}
d((x_1,y_1),(x_2,y_2))=\max\left\{d_{X}(x_1,x_2),d_{Y}(y_1,y_2)\right\}.
\end{equation*}
\end{proposition}

\begin{proof}
Suppose $T$ is restricted sensitive at $x \in X$ with sensitivity constant $\delta>0$ and asymptotic rate $a>0$. For any $\varepsilon>0$, suppose $n \geq $ is such that
\[\mu\left\{x'\in B^{d_X}_\varepsilon(x):d_{X}(T^{n}(x),T^{n}(x'))>\delta\right\}>0.\]
\noindent Then for any $y \in Y$,
\begin{align*}
&(\mu\times\nu)\left\{(x',y')\in B^{d}_{\varepsilon}(x,y): d((T\times S)^{n}(x,y),(T\times S)^{n}(x',y'))>\delta\right\}\\
& \geq \nu (B^{d_Y}_{\varepsilon}(y))\mu\left\{x'\in B^{d_X}_{\varepsilon}(x):d_{X}(T^{n}(x),T^{n}(x'))>\delta\right\}>0.
\end{align*}
As $\nu$ is a probability measure, $(\mu\times\nu)B^{d}_{\varepsilon}(x,y) \leq  \mu B^{d_X}_{\varepsilon}(x)$ for any $\varepsilon>0$, so $-a \log(\mu\times\nu)B^{d}_{\varepsilon}(x,y) \geq  -a \log\mu B^{d_X}_{\varepsilon}(x)$. Hence $T \times S$ is restricted sensitive at $(x,y)$ with the same sensitivity constant $\delta$ and asymptotic rate $a$. This holds for a.e. $x \in X$ and all $y \in Y$, so $T \times S$ is restricted sensitive.
\end{proof}

\begin{corollary}\label{typeIII}
Let $T$ be a measure-preserving mixing (or mildly mixing) transformation on a probability space $(X,\mu)$ with metric $d_X$, and suppose that $T$ is restricted sensitive. Let $S$ be a type III nonsingular conservative ergodic invertible transformation on a probability space $(Y,\nu)$ with metric $d_Y$. Then $T \times S$ on $X \times Y$ with the metric $d$ given in Proposition \ref{prod} is a type III conservative ergodic transformation that is restricted sensitive.
\end{corollary}

\begin{proof} By Proposition 5.4 and Theorem 5.2 of \cite{H-S97} the transformation
$T\times S$ is conservative ergodic and type III. Proposition~\ref{prod} implies that $T\times S$ is
restricted sensitive.
\end{proof}

\bibliographystyle{amsalpha}
\bibliography{RestSen_Bib}

\newcommand{\etalchar}[1]{$^{#1}$}
\providecommand{\bysame}{\leavevmode\hbox to3em{\hrulefill}\thinspace}
\providecommand{\MR}{\relax\ifhmode\unskip\space\fi MR }
\providecommand{\MRhref}[2]{%
  \href{http://www.ams.org/mathscinet-getitem?mr=#1}{#2}
}
\providecommand{\href}[2]{#2}
\begin{thebibliography}{BBC{\etalchar{+}}92}

\bibitem[AAB96]{A-B96}
Ethan Akin, Joseph Auslander, and Kenneth Berg, \emph{When is a transitive map
  chaotic?}, Convergence in ergodic theory and probability (Columbus, OH,
  1993), Ohio State Univ. Math. Res. Inst. Publ., vol.~5, de Gruyter, Berlin,
  1996, pp.~25--40. \MR{MR1412595 (97i:58106)}

\bibitem[ABC02]{ABC02}
Christophe Abraham, G{\'e}rard Biau, and Beno{\^{\i}}t Cadre, \emph{Chaotic
  properties of mappings on a probability space}, J. Math. Anal. Appl.
  \textbf{266} (2002), no.~2, 420--431. \MR{MR1880515 (2002j:37009)}

\bibitem[AFS01]{AFS01}
Terrence Adams, Nathaniel Friedman, and Cesar~E. Silva, \emph{Rank-one power
  weakly mixing non-singular transformations}, Ergodic Theory Dynam. Systems
  \textbf{21} (2001), no.~5, 1321--1332. \MR{MR1855834 (2002h:37013)}

\bibitem[BBC{\etalchar{+}}92]{BBCDS}
J.~Banks, J.~Brooks, G.~Cairns, G.~Davis, and P.~Stacey, \emph{On {D}evaney's
  definition of chaos}, Amer. Math. Monthly \textbf{99} (1992), no.~4,
  332--334. \MR{MR1157223 (93d:54059)}

\bibitem[CJ05]{Cadre05}
Beno{\^{\i}}t Cadre and Pierre Jacob, \emph{On pairwise sensitivity}, J. Math.
  Anal. Appl. \textbf{309} (2005), no.~1, 375--382. \MR{MR2154050
  (2006b:28026)}

\bibitem[CS04]{CS04}
Darren Creutz and Cesar~E. Silva, \emph{Mixing on a class of rank-one
  transformations}, Ergodic Theory Dynam. Systems \textbf{24} (2004), no.~2,
  407--440. \MR{MR2054050 (2005a:37011)}

\bibitem[DS09]{DS09}
Alexandre~I. Danilenko and Cesar~E. Silva, \emph{Ergodic theory: Nonsingular
  transformations}, Encyclopedia of Complexity and System Science, vol. Part 5,
  Springer, 2009, pp.~3055--3083.

\bibitem[GII{\etalchar{+}}08]{G-S08}
Ilya Grigoriev, Nathaniel Ince, C\u{a}t\u{a}lin~M. Iordan, Amos Lubin, and
  Cesar~E. Silva, \emph{On $\mu$-compatible metrics and measurable
  sensitivity}, preprint. (2008).

\bibitem[GW93]{GW93}
Eli Glasner and Benjamin Weiss, \emph{Sensitive dependence on initial
  conditions}, Nonlinearity \textbf{6} (1993), no.~6, 1067--1075. \MR{MR1251259
  (94j:58109)}

\bibitem[HLYar]{Huang}
Wen Huang, Ping Lu, and Xiangdong Ye, \emph{Measure-theoretical sensitivity and
  equicontinuity}, Israel J. Math. (to appear), 1--51.

\bibitem[HS98]{H-S97}
Jane Hawkins and Cesar~E. Silva, \emph{Characterizing mildly mixing actions by
  orbit equivalence of products}, New York J. Math. \textbf{3A} (1997/98),
  no.~Proceedings of the New York Journal of Mathematics Conference, June
  9--13, 1997, 99--115 (electronic). \MR{MR1611117 (99b:28019)}

\bibitem[HYW04]{He04}
Lianfa He, Xinhua Yan, and Lingshu Wang, \emph{Weak-mixing implies sensitive
  dependence}, J. Math. Anal. Appl. \textbf{299} (2004), no.~1, 300--304.
  \MR{MR2091290 (2005i:37005)}

\bibitem[JKL{\etalchar{+}}08]{J-S08}
Jennifer James, Thomas Koberda, Kathryn Lindsey, Cesar~E. Silva, and Peter
  Speh, \emph{Measurable sensitivity}, Proc. Amer. Math. Soc. \textbf{136}
  (2008), no.~10, 3549--3559.

\bibitem[Moo07]{Moot07}
T.~K.~Subrahmonian Moothathu, \emph{Stronger forms of sensitivity for dynamical
  systems}, Nonlinearity \textbf{20} (2007), no.~9, 2115--2126. \MR{MR2351026}

\bibitem[Orn60]{Or60}
Donald~S. Ornstein, \emph{On invariant measures}, Bull. Amer. Math. Soc.
  \textbf{66} (1960), 297--300. \MR{MR0146350 (26 \#3872)}

\bibitem[Orn72]{Or70}
\bysame, \emph{On the root problem in ergodic theory}, Proceedings of the Sixth
  Berkeley Symposium on Mathematical Statistics and Probability (Univ.
  California, Berkeley, Calif., 1970/1971), Vol. II: Probability theory
  (Berkeley, Calif.), Univ. California Press, 1972, pp.~347--356. \MR{MR0399415
  (53 \#3259)}

\bibitem[Sil08]{Si08}
C.~E. Silva, \emph{Invitation to ergodic theory}, Student Mathematical Library,
  vol.~42, American Mathematical Society, Providence, RI, 2008. \MR{MR2371216}

\bibitem[Sri98]{srivastava1998cbs}
S.M. Srivastava, \emph{{A Course on Borel Sets}}, Springer, 1998.

\end{thebibliography}

\appendix
\section{Measurability of the Minimum Asymptotic Rate}\label{appendix}
The proof of Proposition \ref{ameasurable} requires the following lemmas:

\begin{lemma}\label{muIsLSC}
Let $(X,d)$ be a metric space with a probability measure $\mu$. For all $r>0$, the function $f(x)=\mu B_r(x)$ is lower semi-continuous.
\end{lemma}

\begin{proof}
Note that for any $x \in X$,
\[\lim_{n \to \infty} \mu B_{r-\frac{1}{n}}(x)=\mu \left(\bigcup_{n=1}^\infty B_{r-\frac{1}{n}}(x)\right)=\mu B_r(x).\]
Hence for any $\varepsilon>0$, there exists $n$ such that $\mu B_r(x)-\mu B_{r-\frac{1}{n}}(x)<\varepsilon$. For any $y \in B_{\frac{1}{n}}(x)$, $B_{r-\frac{1}{n}}(x) \subset B_r(y)$, so $\mu B_{r-\frac{1}{n}}(x)\leq\mu B_r(y)$ and thus $\mu B_r(x)-\mu B_r(y)<\varepsilon$. Hence $f$ is lower semi-continuous.
\end{proof}

\begin{lemma}\label{Anepsmeas}
Suppose $(X,d)$ is a separable metric space, $\mu$ is a Borel probability measure, and $T$ is a transformation on $X$. Let $\delta>0$, $\varepsilon>0$, and $n \geq 0$, and let $A_{n,\varepsilon,\delta}=\{x \in X : \mu\{y \in B_\varepsilon(x) : d(T^nx,T^ny)>\delta\}>0\}$. Then $A_{n,\varepsilon,\delta}$ is measurable.
\end{lemma}

\begin{proof}
By Lusin's Theorem, for each integer $k>0$, there exists a closed set $F_k \subset X$ such that $\mu(X \setminus F_k)<\frac{1}{k}$ and $T|_{F_k}:F_k \to X$ is continuous. Let $E_k=X \setminus \bigcup_{i=0}^{n-1} T^{-i}(X \setminus F_k)$, so $T^n$ is continuous on $E_k$, and define $f_k:E_k \to [0,1]$ as $f_k(x)=\mu\{y \in B_\varepsilon(x) \cap E_k : d(T^nx,T^ny) \geq \delta\}$.

We claim that $f_k$ is lower-semicontinuous. Suppose $f_k(x)>0$ for some $x \in E_k$, and let $c>0$. There exists $\eta>0$ such that $\mu\{y \in B_\varepsilon(x) \cap E_k : d(T^nx,T^ny) \geq \delta+\eta\}>f_k(x)-\frac{c}{2}$. Since $T|_{E_k}$ is continuous, there exists $r>0$ such that $d(T^nx,T^nx')<\eta$ for all $x' \in B_r(x) \cap E_k$. We may choose $r$ sufficiently small so that $\mu B_{\varepsilon-r}(x)>\mu B_\varepsilon(x)-\frac{c}{2}$. Then for all $x' \in B_r(x) \cap E_k$, $\mu (B_\varepsilon(x) \cap B_\varepsilon(x')) \geq \mu B_{\varepsilon-r}(x)>\mu B_\varepsilon(x)-\frac{c}{2}$. Thus
\begin{align*}
f_k(x')&=\mu\{y \in B_\varepsilon(x') \cap E_k : d(T^nx',T^ny) \geq \delta\}\\
&\geq \mu\{y \in B_\varepsilon(x') \cap E_k : d(T^nx,T^ny) \geq \delta+\eta\}\\
&\geq \mu\{y \in B_\varepsilon(x') \cap B_\varepsilon(x) \cap E_k : d(T^nx,T^ny) \geq \delta+\eta\}\\
&\geq \mu\{y \in B_\varepsilon(x) \cap E_k : d(T^nx,T^ny) \geq \delta+\eta\}-\mu(B_\varepsilon(x) \cap E_k \setminus B_\varepsilon(x'))\\
&>f_{k}(x)-\frac{c}{2}-\frac{c}{2}\\
&=f_{k}(x)-c.
\end{align*}
Hence $f_k$ is lower-semicontinuous on $E_k$, so $f_k$ is measurable since $\mu$ is a Borel measure. Then $f_k^{-1}(0,1]=\{x \in X : \mu\{y \in B_\varepsilon(x) \cap E_k : d(T^nx,T^ny) \geq \delta\} \geq 0\}$ is measurable.

Let $S_{n,\varepsilon,\delta}=\{x \in X: \mu\{y \in B_\varepsilon(x) : d(T^nx,T^ny) \geq \delta\}>0\}$. We note that $f_k^{-1}(0,1] \subset S_{n,\varepsilon,\delta}$ for all $k$. Since $\mu E_k^{c} < \frac{n}{k}$, for each $x \in S_{n,\varepsilon,\delta}$, if there exists $k$ such that $\mu\{y \in B_\varepsilon(x) : d(T^nx,T^ny) \geq \delta\}>\frac{n}{k}$, then $f_k(x)>0$. Hence $x \in f_k^{-1}(0,1]$ for some $k$, so $S_{n,\varepsilon,\delta}=\bigcup_{k=1}^\infty f_k^{-1}(0,1]$, which is measurable. Finally, $A_{n,\varepsilon,\delta}=\bigcup_{\delta'<\delta} S_{n,\varepsilon,\delta'}$, which is measurable.
\end{proof}

\begin{proof}[Proof of Proposition ~\ref{ameasurable}]
For any $\varepsilon>0$ and $\delta>0$, define $n_{\varepsilon,\delta}:X \to \mathbb{Z}$ as $n_{\varepsilon,\delta}(x)=\min \{n \geq 0 : \mu \{y \in B_\varepsilon(x) : d(T^nx,T^ny)>\delta\}>0\}$. For any $n \geq 0$, let $A_{n,\varepsilon,\delta}=\{x \in X : \mu\{y \in B_\varepsilon(x) : d(T^nx,T^ny)>\delta\}>0\}$. Define $S_{n,\varepsilon,\delta}={A_{n,\varepsilon,\delta} \setminus \bigcup_{i=0}^{n-1} A_{i,\varepsilon,\delta}}$. Then $n_{\varepsilon,\delta}(x)=\sum_{n=0}^\infty n \chi_{S_{n,\varepsilon,\delta}}(x)$. By Lemma~\ref{Anepsmeas}, $A_{n,\varepsilon,\delta}$ is measurable for all $n$, so $n_{\varepsilon,\delta}$ is measurable. For any $\delta>0$, consider
\begin{align*}
\bar{a}_\delta(x) = \inf\{&a>0 : \forall \varepsilon>0,\,\exists n \leq -a \log \mu B_\varepsilon(x) \text{ s.t. }\\
\nonumber &\mu\{y \in B_\varepsilon(x) : d(T^nx,T^ny)>\delta\}>0\},
\end{align*}
as an extended real-valued function with $\bar{a}_\delta(x)=\infty$ if the set of such $a>0$ is empty (i.e. if $\delta$ is not a sensitivity constant for $x \in X$). Then $\bar{a}_\delta(x)=\sup_{\varepsilon>0} \frac{n_{\varepsilon,\delta}(x)}{-\log \mu B_\varepsilon(x)}$. By Lemma~\ref{muIsLSC}, $-\log \mu B_\varepsilon(x)$ is upper semi-continuous and thus measurable since $\mu$ is a Borel measure, so $\frac{n_{\varepsilon,\delta}(x)}{-\log \mu B_\varepsilon(x)}$ is measurable.

Let us fix $x$ and $\delta$ and consider $n_{\varepsilon,\delta}(x)$ and $-\log \mu B_\varepsilon(x)$ as monotonically decreasing functions of $\varepsilon$. We note that both functions are left-continuous because $\mu B_\varepsilon(x)$ is left-continuous in $\varepsilon$. For each $n$, let $\varepsilon_{\delta,n}=\max \{\varepsilon>0 : n_{\varepsilon,\delta}(x)=n\}$. Then $n_{\varepsilon,\delta}(x)$ as a function of $\varepsilon$ is constant on each of the intervals $(\varepsilon_{n,\delta},\varepsilon_{n+1,\delta}]$ for all $n$, so $\frac{n_{\varepsilon,\delta}(x)}{-\log \mu B_\varepsilon(x)}$ as a function of $\varepsilon$ is monotonically increasing on these intervals. Hence $\bar{a}_\delta(x)=\sup_{n \geq 0} \lim_{\varepsilon \to \varepsilon_{n,\delta}^-} \left(\frac{n_{\varepsilon,\delta}(x)}{-\log \mu B_\varepsilon(x)}\right)$, which is measurable. Finally, we note that for any fixed $x$, $\bar{a}_\delta(x)$ is monotonically increasing in $\delta$, so $a_T^*=\inf_{\delta>0} \bar{a}_\delta=\lim_{\delta \to 0} \bar{a}_\delta$ as extended real-valued functions, which is measurable.
\end{proof}


\end{document}